\documentclass{article}

\RequirePackage[english]{babel}
\RequirePackage{a4wide}
\usepackage{float}
\usepackage{imakeidx}
\makeindex

\RequirePackage{amsmath,amsfonts,amssymb,amsthm,mathtools,bm,siunitx}
\RequirePackage{graphicx,xcolor}
\graphicspath{{./fig/}}
\usepackage{ulem}

\newcommand{\sib}[1]{[\si{#1}]}

\theoremstyle{plain}

\title{Sensitivity analysis with a 3D mixed-dimensional code for DC geoelectrical investigations of landfills: synthetic tests}

\author{
Lorenzo Panzeri$^1$ \and
Alessio Fumagalli$^{1, *}$ \and
Laura Longoni$^1$ \and
Monica Papini$^1$ \and
Diego Arosio$^2$}

\date{$^1$ Politecnico di Milano\\%
$^2$ Universit\`a degli Studi di Modena e Reggio Emilia\\
$^*$ alessio.fumagalli@polimi.it}

\begin{document}

\maketitle

\begin{abstract}
Electrical resistivity tomography is a suitable technique for non-invasive monitoring of municipal solid waste landfills, but accurate sensitivity analysis is necessary to evaluate the effectiveness and reliability of geoelectrical investigations and to properly design data acquisition.
Commonly, a thin high-resistivity membrane in placed underneath the waste to prevent leachate leakage.
In the construction of a numerical framework for sensitivity computation, taking into account the actual dimensions of the electrodes and, in particular, of the membrane, can lead to extremely high computational costs.
In this work, we present a novel approach for numerically computing sensitivity effectively by adopting a mixed-dimensional framework, where the membrane is approximated as a 2D object and the electrodes as 1D objects. The code is first validated against analytical expressions for simple 4-electrode arrays and a homogeneous medium. It is then tested in simplified landfill models, where a 2D box-shaped liner separates the landfill body from the surrounding media, and 48 electrodes are used.
The results show that electrodes arranged linearly along both sides of the perimeter edges of the box-shaped liner are promising for detecting liner damage, with sensitivity increasing by 2-3 orders of magnitude, even for damage as small as one-sixth of the electrode spacing in diameter. Good results are also obtained when simulating an electrical connection between the landfill and the surrounding media that is not due to liner damage.
The next steps involve evaluating the minimum number of configurations needed to achieve suitable sensitivity with a manageable field effort and validating the modeling results with downscaled laboratory tests.

\end{abstract}

\section{Introduction}

Nowadays, the creation of new municipal solid waste landfills (MSWLFs) is discouraged in favour of a circular economy, with minimum and maximum rates for recycling and landfilling, respectively, set by regulatory agencies, especially in developed countries (\cite{EC2015}). As a consequence, landfilling of municipal solid waste is declining, also thanks to the advancements made in recycling, composting, incineration and energy recovery technologies (\cite{Nanda2021}). Nevertheless, it is estimated that there are hundreds of thousands of active and closed landfills in the European Union (\cite{Jones2013,Wagner2015}), as well as in the USA and Canada (\cite{Giroux2014,Peters2016}).

In recent MSWLFs, besides facilities to collect leachate and gas produced by waste degradation, one or more high-density polyethylene (HDPE) liners are placed underneath the waste body to prevent leachate leakage in the subsurface and potential groundwater contamination (\cite{Nanda2021}). Monitoring the conditions of the liners and of the materials beneath the waste is crucial in landfill management to avoid serious environmental concerns.

DC geoelectrical surveys are well suited to investigate the electrical resistivity of the media and can be effective in MSWLF monitoring as plastic liner is highly resistive and leachate in very conductive because of high salt content. Both Electrical Resistivity Tomography (ERT) and Induced Polarization (IP) have been employed for non-invasive assessment and monitoring of landfills in field (\cite{Dahlin2010,DeCarlo2013,Tsourlos2014,DEDONNO2017302}), laboratory and numerical studies (\cite{Frangos1997,Binley2003,Ling2019,Aguzzoli2020,panzeri2023validation}).

ERT investigations are typically conducted using multi-electrode resistivity meters capable of automatically performing measurements with several quadrupoles. Each quadrupole consists of four electrodes: two are used to inject electrical current into the subsurface, while the other two measure the voltage difference. Besides the signal-to-noise ratio, the geometry of the quadrupole controls the investigation depth and spatial resolution (\cite{Loke2022}) and can be modified to tailor the sensitivity of geoelectrical measurements to specific subsurface features or target depths. Integrating data from multiple quadrupoles helps to improve the accuracy of the resulting resistivity images and thus offers enhanced imaging capabilities of subsurface structures.

However, considering the large size of the landfills and the fact that the electrodes are generally placed around or above the waste body, the depth of investigation and the spatial resolution of the geoelectrical technique can be very limited where they are most needed, namely below the waste mass. It is well known that the sensitivity of the resistivity investigation rapidly decays away from the deployed electrodes (\cite{Loke2022}).

    In this work, we consider the sensitivity of ERT within the framework of a 3D mixed-dimensional numerical model developed for MSWLF monitoring to reduce computational costs and address the ill-posedness of the geoelectrical inverse problem (\cite{Fumagalli2023}). The paper is organized as follows: In Section 2, we review the mixed-dimensional mathematical model, where the electrodes and the HDPE liner are modeled as 1D and 2D elements, respectively. In Section 3, we first introduce the classical equations of geoelectrical sensitivity and then extend them to the newly developed mixed-dimensional framework. Section 4 describes some practical implementation aspects of the code used for computing sensitivity. Section 5 presents numerical experiments, where sensitivity is first computed for validation purposes and then in the presence of simplified landfill models. The results of the numerical experiments are discussed in Section 6. Finally, conclusions are drawn in Section 7.

    The ultimate goal of our work is to provide a tool for computing DC sensitivity and evaluating the effectiveness of geoelectrical investigations of MSWLFs. We aim to contribute to the advancement of non-destructive techniques for assessing the condition of HDPE liners and potential contaminant plumes at MSWLF sites, with implications for enhanced environmental protection and sustainable waste management practices.

\section{The mathematical  model}\label{sec:mathmolde}

We adopt the mixed-dimensional model proposed by \cite{Fumagalli2023}. In this section, for completeness, we recall the main points and the equations that are useful for the derivation of the sensitivity presented in Section \ref{sec:sensitivity}.

We indicate the domain, where the geoelectrical equations will be applied, as $\Upsilon\subset \mathbb{R}^3$, with boundary $\partial \Upsilon$ and outward unit normal $\bm{n}_\partial$. We consider the following main variables: $\bm{J}:\Upsilon \rightarrow \mathbb{R}^3$, the current density in \sib{\ampere\per\square\meter}, and
$\varphi: \Upsilon \rightarrow \mathbb{R}$, the electric potential in \sib{\volt}. We indicate with $\bm{x} = (x_0, x_1, x_2)$ a generic point.

The constitutive relation between the current density $\bm{J}$ and the electric potential $\varphi$ is given by Ohm's law as follows
\begin{gather}\label{eq:current_potential}
    \rho \bm{J} + \nabla \varphi = \bm{0}.
\end{gather}
where $\rho$ is the electrical resistivity in \sib{\ohm\meter}. Gauss' law, or charge conservation equation for $\bm{J}$, can be expressed as
\begin{gather}\label{eq:continuity}
    \nabla \cdot \bm{J} = q,
\end{gather}
with $q:\Upsilon \rightarrow \mathbb{R}$ being the source of volumetric charge
density in \sib{\ampere\per\cubic\meter}.

In the considered examples, $\Upsilon$ is the domain defined as
$\Upsilon = \mathbb{R}\times \mathbb{R}\times \mathbb{T}$, where $\mathbb{T}$ is the hypograph of a smooth function $x_2 = t(x_0, x_1)$ representing the elevation associated to a digital
terrain model, so $\mathbb{T} = \text{hyp}(t)$. For a flat terrain we have $x_2 = \overline{x}_2$, with
$\overline{x}_2 \in \mathbb{R}$, and thus we get
$\mathbb{T} = (-\infty, \overline{x}_2)$ and so $\Upsilon =  \mathbb{R}
    \times \mathbb{R} \times (-\infty, \overline{x}_2)$. However, this might be unpractical since unbounded domains are complicated to be considered in numerical simulations. We know that if $q$ has limited support, then it can be shown that
\begin{gather}\label{eq:at_inf}
    \lim_{ \Vert \bm{x} \Vert \rightarrow \infty} \varphi(\bm{x}) = 0
\end{gather}
and, with a reasonable approximation, we restrict $\Upsilon$ to be a bounded domain and impose that at its lateral and bottom boundaries the electric potential is null. On the top boundary, we have air that, for our purposes, does not interact with $\Upsilon$ and we can assume $\bm{J}\cdot\bm{n}_\partial = 0$ on the top. We decompose $\partial \Upsilon$ into two non-overlapping parts     $\partial^{t} \Upsilon$ and $\partial^o \Upsilon$, representing the top boundary (i.e., the digital terrain model) and the lateral and bottom boundaries, respectively.

By combining \eqref{eq:current_potential} and \eqref{eq:continuity} the direct current problem can be written as:
find $(\bm{J}, \varphi)$ in $\Upsilon$ such that
\begin{gather}\label{eq:the_system}
    \begin{aligned}
         & \begin{aligned}
             & \rho \bm{J} + \nabla \varphi = \bm{0} \\
             & \nabla \cdot \bm{J} = q
        \end{aligned}
         &
         &                           & \text{in } \Upsilon            \qquad
        \begin{aligned}
             & \varphi = 0
             &                                  & \text{on } \partial^o \Upsilon \\
             & \bm{J} \cdot \bm{n}_\partial = 0
             &                                  & \text{on } \partial^t \Upsilon
        \end{aligned}
    \end{aligned}
\end{gather}
The previous problem can be written in its primal formulation, setting $\sigma = \rho^{-1}$ the electric
conductivity in $\sib{\siemens\per\meter}$ and, by substituting Ohm's law into we obtain the conservation equation, we get
\begin{gather*}
    \begin{aligned}
         & \nabla \cdot (-\sigma \nabla \varphi) = q
         &                                           & \text{in } \Upsilon \qquad
        \begin{aligned}
             & \varphi = 0
             &                                                  & \text{on } \partial^o \Upsilon \\
             & -\sigma \nabla \varphi \cdot \bm{n}_\partial = 0
             &                                                  & \text{on } \partial^t \Upsilon
        \end{aligned}
    \end{aligned}
\end{gather*}

It is possible to compute the solution $\varphi$ at a point $\bm{y}$, by using the $\sigma$-weighted
Green's
function $g^{\bm{y}}$, in \sib{\ohm}, related to the point $\bm{y}$. The Green function
solves the following problem and is such that
\begin{gather}\label{eq:green}
    \nabla \cdot(-\sigma \nabla g^{\bm{y}}) = \delta^{\bm{y}}
    \quad \Rightarrow \quad
    \varphi (\bm{y}) = \int_\Upsilon {q}  g^{\bm{y}}
\end{gather}
where $\delta^{\bm{y}}$ is a point delta centered in $\bm{y}$ measured in $\sib{\per\cubic\meter}$.
The green function is a fundamental tool in solving problems
related to potential theory and electrostatics. It represents
the potential due to a unit point charge located at the point $\bm{y} \in \Upsilon$ and, if
$\Upsilon$ is unbounded, is such that
\begin{gather*}
    \lim_{\Vert \bm{x} \Vert \rightarrow \infty} g^{\bm{y}}(\bm{x})  = 0.
\end{gather*}
For a homogeneous three-dimensional half space, electrodes at the surface and unitary injected current, the analytical expression for $g^{\bm{y}}$ is well known
\begin{gather*}
    g^{\bm{y}}(\bm{x}) =\frac{\rho}{2\pi \Vert \bm{x} - \bm{y}\Vert},
\end{gather*}
On the contrary, except a few cases (\cite{Sheriff1990}), the Green function has to be approximated when the half space is not homogeneous.

In MSWLFs geoelectrical investigations, the electrodes used to inject current and measure potential difference can generally be approximated as cylinders whose radius is much smaller than their height. Moreover, the low-permeability HDPE liner (hereafter referred to as $\lambda$) used to seal the waste body has a thickness ($\epsilon$) of a few millimeters and covers an area of thousands of square meters. Accordingly, it is numerically impractical to represent these objects with a 3D computational grid, so we follow the strategy proposed by \cite{Fumagalli2023} and approximate them as lower-dimensional objects: the electrodes are represented as one-dimensional and the liner as a two-dimensional object. The grid now does not resolve the radius of the electrodes nor the thickness of the liner and can be coarsened for computational efficiency.

This approach, however, requires to consider dedicated mathematical models. More in detail, for a single electrode $\gamma$, the equations to compute $\bm{J}_\gamma$ and $\varphi_\gamma$, now expressed in $\sib{\ampere}$ and $\sib{\volt}$, are summarized as: find $(j_\gamma, \bm{J}_\gamma, \varphi_\gamma)$ in $\gamma$ such that
\begin{gather}\label{eq:the_system_electrode}
    \begin{aligned}
         & \begin{aligned}
             & \rho_\gamma \bm{J}_\gamma + \pi r^2 \nabla \varphi_\gamma = \bm{0} \\
             & \nabla \cdot \bm{J}_\gamma - j_\gamma= 0
        \end{aligned}
         &
         &                            & \text{in } \gamma            \qquad
        \begin{aligned}
             & \bm{J}_\gamma \cdot \bm{n}_\partial = 0
             &                                                            & \text{on } \partial^o \gamma                     \\
             & \bm{J}_\gamma \cdot \bm{n}_\partial = \pi r^2 \overline{J}
             &                                                            & \text{on } \partial^t \gamma                     \\
             & \rho_\gamma j_\gamma + \varphi_\gamma - \varphi = 0        &                              & \text{on } \gamma
        \end{aligned}
    \end{aligned}
\end{gather}
where $j_\gamma$ in $\sib{\ampere\per\meter}$ is the current density exchanged between the electrode and the surrounding domain $\Upsilon$ and $\rho_\gamma$ in $\sib{\ohm\,\meter}$ is the resistivity of the electrode. The top boundary of the electrode $\partial^t \gamma$ is in contact with the top boundary of the domain $\partial^t \Upsilon$, while $\partial^o \gamma$ is the portion of the electrode boundary immersed into $\Upsilon$.

Having assumed the 1D representation of a cylindrical electrode $\gamma$, in \eqref{eq:the_system_electrode} it can be seen that Ohm's low in \eqref{eq:current_potential} is now averaged over each cross-section of radius $r$.

The source term $q$ introduced in \eqref{eq:the_system} is now given by $q =- j_\gamma \delta_\gamma$, where $\delta_\gamma$ in $\sib{\per\square\meter}$ is a linear delta function distributed along $\gamma$. This
couples problems \eqref{eq:the_system} and \eqref{eq:the_system_electrode}, being $j_\gamma$ the coupling variable. Moreover, the boundary
condition in \eqref{eq:the_system_electrode} on $\partial^o \gamma$ is the tip-condition
(no-current flow) at the point immersed in $\Upsilon$, while on $\partial^t \gamma$  we impose the injected current $\overline{J}$, measured in \sib{\ampere\per\square\meter}.

Also for the liner $\lambda$, represented as a two-dimensional object, we need a dedicated model to compute $\bm{J}_\lambda$ in $\sib{\ampere\per\meter}$ and $\varphi_\lambda$ in $\sib{\volt}$. Thus, we consider the following problem summarized as: find $(j_\lambda, \bm{J}_\lambda, \varphi_\lambda)$ in
$\lambda$ such that
\begin{gather}\label{eq:the_system_liner}
    \begin{aligned}
         & \begin{aligned}
             & \rho_\lambda \bm{J}_\lambda + \epsilon \nabla \varphi_\lambda = \bm{0} \\
             & \nabla \cdot \bm{J}_\lambda - j_\lambda= 0
        \end{aligned}
         &
         &                            & \text{in } \lambda            \qquad
        \begin{aligned}
             & \bm{J}_\lambda \cdot \bm{n}_\partial = 0
             &                                                                  & \text{on } \partial \lambda                      \\
             & \epsilon \rho_\lambda  j_\lambda + \varphi_\lambda - \varphi = 0 &                             & \text{on } \lambda
        \end{aligned}
    \end{aligned}
\end{gather}
where $j_\lambda$ in $\sib{\ampere\per\square\meter}$ and $\rho_\lambda$ in $\sib{\ohm\,\meter}$ are the current density exchanged between the liner and its surrounding and the resistivity of the membrane, respectively.

In \eqref{eq:the_system_liner}, it has been assumed that the liner has thickness much smaller than its other dimensions, and so the equations had been integrated along every cross section of thickness $\epsilon$. It is also assumed that the boundary of the liner $\partial \lambda$ does not exchange any current density with the surrounding.

\section{Sensitivity} \label{sec:sensitivity}

In geoelectrical sounding, the sensitivity is just a measure of how much the electric potential (i.e., the observed data) changes in response to a change in resistivity (i.e., the model parameters) within the considered domain. In other words, the sensitivity is fully described by the Jacobian matrix containing the first derivatives of the potential with respect to the model resistivities. Here, the sensitivity analysis is performed through the use of Green's functions, as mentioned in \cite{mcgillivray1990methods} and \cite{park1991inversion}. In the next two subsections, we first recall the computation of sensitivity in a standard equi-dimensional model and then we apply the same approach to our mixed-dimensional model.

\subsection{Equi-dimensional sensitivity}

Let us consider \eqref{eq:the_system} with $q = I \delta^{\bm{x}}$, which represents an electrode $\gamma$ placed at position $\bm{x}$ and with $I$ the current exchanged. Accordingly, in
\eqref{eq:the_system_electrode} we have $I = j_\gamma$.

Following the discussion in \cite{mcgillivray1990methods}, we assume that
$\Upsilon$ can be divided into $N$ non-overlapping subdomains $\omega_i$, that are the cells of the model and will form the computational grid. We can express the resistivity, and similarly the electric conductivity, and its derivative as
\begin{gather*}
    \rho(\bm{x}) = \sum_{k=1}^N \rho_k \chi_k (\bm{x}) \quad \text{and}
    \quad {\partial_{\rho_i}}\rho (\bm{x}) = \chi_i (\bm{x})
\end{gather*}
where $\partial_{\rho_i} \star= \frac{\partial \star}{\partial \rho_i}$ for brevity and $\chi_i$ is the dimensionless indication function of the $i$-th subdomain defined as
\begin{gather*}
    \chi_i (\bm{x}) =
    \begin{cases}
        1 & \text{if } \bm{x} \in \omega_i \\
        0 & \text{otherwise}
    \end{cases}
\end{gather*}
Since both variables depend on the resistivity of the system $\bm{J} = \bm{J}(\rho)$ and
$\varphi = \varphi(\rho)$,
we define the sensitivity associated to the resistivity in the subdomain $\omega_i$ as
\begin{gather}\label{eq:sens_defi}
    \bm{J}_i = \partial_{\rho_i} \bm{J}
    \quad \text{and} \quad
    \varphi_i = \partial_{\rho_i} \varphi
\end{gather}
with dimension \sib{\ampere\siemens\per\cubic\meter} and \sib{\volt\siemens\per\meter}=\sib{\ampere\per\meter},
respectively.
We now differentiate the first equation of \eqref{eq:the_system} to get the following expression
\begin{gather*}
    \partial_{\rho_i} (\rho\bm{J} + \nabla \varphi) =   \chi_i \bm{J} + \rho \bm{J}_i
    + \nabla (\partial_{\rho_i} \varphi) =
    \chi_i \bm{J} + \rho \bm{J}_i + \nabla \varphi_i = \bm{0}
\end{gather*}
while the second equation of \eqref{eq:the_system} becomes
\begin{gather*}
    \partial_{\rho_i} \nabla \cdot \bm{J} = \partial_{\rho_i} (I \delta^{\bm{x}})
    \quad \Rightarrow \quad
    \nabla \cdot \bm{J}_i = 0
\end{gather*}
Thus the system for computing the sensitivity of the current density and the electric potential
is given by: find $(\bm{J}_i, \varphi_i)$ such that
\begin{gather}\label{eq:sens}
    \begin{cases}
        \rho \bm{J}_i + \nabla \varphi_i = - \chi_i \bm{J} \\
        \nabla \cdot \bm{J}_i = 0
    \end{cases}
    \text{in } \Upsilon
\end{gather}
The above equations correspond to the equation (70a) in \cite{mcgillivray1990methods} and can also be written in the primal formulation (by substituting the first equation into the second)
\begin{gather*}
    \nabla \cdot (-\sigma \nabla \varphi_i) = \nabla \cdot ( \chi_i \sigma \bm{J})
\end{gather*}
We suppose that a second electrode is used to measure the electric potential at point $\bm{y}$, thus we need to compute the sensitivity not in all the domain but only
in correspondence to that electrode. We get that, by considering \eqref{eq:green} now applied
to $\varphi_i$, the expression is
\begin{gather}\label{eq:sensi}
    \varphi_i (\bm{y}) = \int_\Upsilon \nabla \cdot (\chi_i \sigma \bm{J}) g^{\bm{y}}
    = -\int_\Upsilon \chi_i \sigma \bm{J}  \cdot \nabla g^{\bm{y}} =
    \sigma_i^2\int_{\omega_i} \nabla \varphi  \cdot \nabla g^{\bm{y}} =
    I \sigma_i^2\int_{\omega_i}  \nabla g^{\bm{x}} \cdot \nabla g^{\bm{y}}
\end{gather}
which is an easy way to compute the sensitivity of the $i$-th subdomain, measured at the
observation point located in $\bm{y}$, since $\overline{J}$ is set in
\eqref{eq:the_system_electrode} and the previous expression depends only on the computation of the Green functions associated with the two electrodes. In other terms, \eqref{eq:sensi} is the relation to compute the sensitivity of a pole-pole configuration (\cite{park1991inversion,loke1995least}).

Now, let us consider the case of a quadrupole, where two electrodes are used to circulate the current density, while the other two measure the difference in electric potential. One of the current electrodes is used to inject current density at point $\bm{x}_0$, while the other one extract the same amount at point $\bm{x}_1$, clearly being $\bm{x}_0 \neq \bm{x}_1$.
The computation of the electric potential requires to set $q =  I (\delta^{\bm{x}_0}  - \delta^{\bm{x}_1})$ in \eqref{eq:the_system}.
The computation of the sensitivity is now associated with the difference of $\varphi$ at the
two measuring electrodes at position $\bm{y}_0$ and $\bm{y}_1$, with $\bm{y}_0 \neq \bm{y}_1$.
We have
\begin{gather*}
    \varphi_i({\bm{y}_0}) - \varphi_i({\bm{y}_1}) =-
    \int_\Upsilon \chi_i \sigma \bm{J} \cdot  \nabla g^{\bm{y}_0} + \int_\Upsilon \chi_i \sigma \bm{J} \cdot \nabla g^{\bm{y}_1}
    = \sigma_i^2\int_{\omega_i} \nabla \varphi \cdot \nabla g^{\bm{y}_0} - \sigma_i^2\int_{\omega_i} \nabla \varphi \cdot \nabla g^{\bm{y}_1}
\end{gather*}
since $\varphi$ is computed with two delta functions, we can use again the $\sigma$-weighted Green functions to
express it as $\varphi = I (g^{\bm{x}_0} - g^{\bm{x}_1})$ and thus the sensitivity can be computed as linear combinations of pole-pole sensitivity functions \eqref{eq:sensi} as
\begin{gather}\label{eq:quadrupole}
    \begin{aligned}
        \varphi_i({\bm{y}_0}) - \varphi_i({\bm{y}_1}) = & +
        I\sigma_i^2 \int_{\omega_i}\nabla g^{\bm{x}_0} \cdot \nabla g^{\bm{y}_0} -
        I\sigma_i^2 \int_{\omega_i} \nabla g^{\bm{x}_1} \cdot \nabla g^{\bm{y}_0}
        \\ &-
        I\sigma_i^2 \int_{\omega_i} \nabla g^{\bm{x}_0} \cdot \nabla g^{\bm{y}_1} +
        I\sigma_i^2 \int_{\omega_i} \nabla g^{\bm{x}_1} \cdot \nabla g^{\bm{y}_1}.
    \end{aligned}
\end{gather}
In the case of multiple electrode configurations, a simple way to obtain the (approximate) global sensitivity is to sum the absolute values of the sensitivities associated with each configuration (\cite{Loke2022}).

\subsection{Mixed-dimensional sensitivity}

Due to presence of the liner, we need to consider the computation of the sensitivity in the mixed-dimensional
framework. As done in the previous section, in
\eqref{eq:the_system} we set $q=I \delta^{\bm{x}}$ and consider also \eqref{eq:the_system_liner}
to describe the liner. In addition to the sensitivity associated to the variables in the 3D domain $\Upsilon$ expressed in \eqref{eq:sens_defi}, we need to consider also the followings associated to the liner
\begin{gather*}
    j_{\lambda,i} = \partial_{\rho_{i}} j_\lambda
    \quad \text{and} \quad
    \bm{J}_{\lambda,i} = \partial_{\rho_{i}} \bm{J}_\lambda
    \quad \text{and} \quad
    \varphi_{\lambda,i} = \partial_{\rho_{i}} \varphi_\lambda
\end{gather*}
where the the variation of $\rho_{i}$ is associated to the
subdomain $\omega_{i}$, which can be part of $\Upsilon$ or $\lambda$. The sensitivities have dimensions $\sib{\ampere\siemens\per\cubic\meter}$,
$\sib{\ampere\siemens\per\square\meter}$ and $\sib{\volt\siemens\per\meter}$=\sib{\ampere\per\meter}, respectively.
Given the partitioning of $\lambda$, the resistivity can be written in the following way
\begin{gather*}
    \rho_\lambda(\bm{x}) = \sum_{k=1}^{N_\lambda} \rho_{\lambda, k} \chi_{\lambda, k} (\bm{x})
\end{gather*}
with $N_\lambda$ is the number of subdomains that compose the liner and $\chi_{\lambda, k}$ is the indication function of the $k$-th subdomain $\omega_{\lambda, k}$.
Let us first suppose that the variation of the resistivity $\rho_i$ is associated with a subdomain $\omega_i$ that
belongs to $\Upsilon$. We have
\begin{gather}\label{eq:sens_lin}
    \begin{cases}
        \rho \bm{J}_i + \nabla \varphi_i = - \chi_i \bm{J} \\
        \nabla \cdot \bm{J}_i = 0
    \end{cases}
    \text{in } \Upsilon
    \qquad
    \begin{cases}
        \rho_\lambda \bm{J}_{\lambda,i} +  \epsilon \nabla \varphi_{\lambda,i} =\bm{0}
        \\
        \nabla \cdot \bm{J}_{\lambda,i} - j_{\lambda,i} = 0
    \end{cases}
    \text{in } \lambda
\end{gather}
while for the interface condition we get
\begin{gather*}
    \partial_{\rho_i}(\epsilon \rho_\lambda  j_\lambda + \varphi_\lambda - \varphi) =
    \epsilon \rho_\lambda \partial_{\rho_i} j_\lambda +
    \partial_{\rho_i} \varphi_\lambda -
    \partial_{\rho_i} \varphi =
    \epsilon \rho_\lambda j_{\lambda,i} +
    \varphi_{\lambda,i} - \varphi_i = 0
\end{gather*}
The above set of equations can also be written in primal form as
\begin{gather*}
    \begin{dcases}
        \nabla \cdot (-\sigma \nabla \varphi_i) = \nabla \cdot ( \chi_i \sigma \bm{J}) & \text{in } \Upsilon \\
        \nabla \cdot (-\epsilon \sigma_\lambda \nabla \varphi_{\lambda, i}) -  j_{\lambda,i} = 0
                                                                                       & \text{in } \lambda  \\
        \epsilon \rho_\lambda j_{\lambda,i} +
        \varphi_{\lambda,i} - \varphi_i = 0
                                                                                       & \text{on } \lambda
    \end{dcases}
\end{gather*}
Let us now extend the concept of the Green function to a mixed-dimensional framework. We call the compound
$(g^{\bm{y}}, g^{\bm{y}}_\lambda, w^{\bm{y}}_\lambda)$, in \sib{\ohm}, \sib{\ohm}, and \sib{\per\square\meter} respectively, the mixed-dimensional $\sigma$-weighted Green function related to a point $\bm{y}$, where a measuring electrode has been placed. The above compound is the solution of the following problem
\begin{gather}\label{eq:greenmd}
    \begin{dcases}
        \nabla \cdot(-\sigma \nabla g^{\bm{y}}) = \delta^{\bm{y}}                                    & \text{in } \Upsilon \\
        \nabla \cdot (-\epsilon \sigma_\lambda \nabla g^{\bm{y}}_\lambda ) - w^{\bm{y}}_\lambda =  0 & \text{in } \lambda  \\
        \epsilon \rho_\lambda w^{\bm{y}}_\lambda + g^{\bm{y}}_\lambda - g^{\bm{y}} = 0               & \text{on }\lambda
    \end{dcases}
\end{gather}
The representation theorem \cite{Salsa2016} can be extended to the mixed-dimensional case and
gives us the possibility to evaluate
the sensitivity of the electric potential. If $q_\lambda$ and $\varpi_{\lambda}$ are the
source terms associated to the equation in $\lambda$ and at the interface between $\lambda$ and its surrounding, we have
\begin{gather*}
    \varphi_i(\bm{y})  =
    \int_\Upsilon q g^{\bm{y}}
    + \int_\lambda q_\lambda g^{\bm{y}}_\lambda
    - \int_\lambda \varpi_{\lambda} w^{\bm{y}}_\lambda
    =
    \int_\Upsilon \nabla \cdot (\chi_i \sigma \bm{J}) g^{\bm{y}}
    =\sigma_i^2 \int_{\omega_i} \nabla \varphi \cdot \nabla g^{\bm{y}}
\end{gather*}
and since $\varphi$ can be expressed in term of $I g^{\bm{x}}$, we obtain the following expression for the sensitivity
\begin{gather}\label{eq:sens_mixed_out}
    \varphi_i(\bm{y})     =I \sigma_i^2 \int_{\omega_i} \nabla g^{\bm{x}} \cdot \nabla g^{\bm{y}}.
\end{gather}

Let us assume now that the variation of the resistivity is associated to a sub-domain of the
liner, named $\omega_{\lambda, i} \subset \lambda$, then the equations for the sensitivity are now given by
\begin{gather}\label{eq:sens_lin_b}
    \begin{cases}
        \rho \bm{J}_i + \nabla \varphi_i = \bm{0} \\
        \nabla \cdot \bm{J}_i = 0
    \end{cases}
    \text{in } \Upsilon
    \qquad
    \begin{cases}
        \rho_\lambda \bm{J}_{\lambda,i} +  \epsilon \nabla \varphi_{\lambda,i} =- \chi_{\lambda, i} \bm{J}_\lambda
        \\
        \nabla \cdot \bm{J}_{\lambda,i} - j_{\lambda,i} = 0
    \end{cases}
    \text{in } \lambda
\end{gather}
while for the interface condition we get
\begin{gather*}
    \partial_{\rho_i}(\epsilon \rho_\lambda  j_\lambda + \varphi_\lambda - \varphi) =
    \epsilon \chi_{\lambda, i}j_\lambda +
    \epsilon \rho_\lambda \partial_{\rho_i} j_\lambda +
    \partial_{\rho_i} \varphi_\lambda -
    \partial_{\rho_i} \varphi =
    \epsilon \chi_{\lambda, i}j_\lambda +
    \epsilon \rho_\lambda j_{\lambda,i} +
    \varphi_{\lambda,i} - \varphi_i = 0
\end{gather*}
The above set of equations can also be written in primal form, in terms only of $(\varphi_i, \varphi_{\lambda,i}, j_{\lambda, i})$, as
\begin{gather*}
    \begin{dcases}
        \nabla \cdot (-\sigma \nabla \varphi_i) = 0 & \text{in } \Upsilon \\
        \nabla \cdot (-\epsilon \sigma_\lambda \nabla \varphi_{\lambda, i}) -  j_{\lambda,i} = \nabla \cdot ( \chi_{\lambda, i} \sigma_{\lambda} \bm{J}_\lambda)
                                                    & \text{in } \lambda  \\
        \epsilon \rho_\lambda j_{\lambda,i} +
        \varphi_{\lambda,i} - \varphi_i =- \epsilon \chi_{\lambda, i}j_\lambda
                                                    & \text{on } \lambda
    \end{dcases}
\end{gather*}
and thus the sensitivity can be computed now with the following expression
\begin{gather*}
    \begin{aligned}
        \varphi_i(\bm{y}) & =
        \int_\Upsilon q g^{\bm{y}}
        + \int_\lambda q_\lambda g^{\bm{y}}_\lambda
        - \int_\lambda \varpi_{\lambda} w^{\bm{y}}_\lambda
        =
        \int_\lambda \nabla \cdot (\chi_{\lambda, i} \sigma_{\lambda} \bm{J}_\lambda ) g^{\bm{y}}_\lambda
        +\int_\lambda\epsilon\chi_{\lambda, i} j_\lambda w^{\bm{y}}_\lambda
        \\
                          & =        - \int_\lambda \chi_{\lambda, i} \sigma_{\lambda} \bm{J}_\lambda \cdot \nabla g^{\bm{y}}_\lambda
        +\epsilon\int_{\omega_{\lambda, i}} j_\lambda w^{\bm{y}}_\lambda
        = \epsilon\sigma_{\lambda,i}^2  \int_{\omega_{\lambda, i}} \nabla \varphi_\lambda \cdot \nabla g^{\bm{y}}_\lambda
        +\epsilon\int_{\omega_{\lambda, i}} j_\lambda w^{\bm{y}}_\lambda
    \end{aligned}
\end{gather*}
Finally, since $\varphi_\lambda$ and $j_\lambda$ can be expressed in term of the Green function at $\bm{x}$, we get
\begin{gather}\label{eq:sens_mixed_int}
    \varphi_i(\bm{y})
    = I\epsilon\sigma_{\lambda,i}^2  \int_{\omega_{\lambda, i}} \nabla g^{\bm{x}}_\lambda \cdot \nabla g^{\bm{y}}_\lambda
    +I \epsilon\int_{\omega_{\lambda, i}} w^{\bm{x}}_\lambda w^{\bm{y}}_\lambda
\end{gather}
When a quadrupole is considered, the previous expressions can be extended in a similar manner to \eqref{eq:quadrupole}.

\section{Numerical model}

In this section we discuss how to numerically approximate the computation of the sensitivity described previously. In particular, in Subsection \ref{subsec:discr} we discuss the numerical scheme adopted and the strategy to couple the grid across the different dimensions.
Subsection \ref{subsec:imple} is devoted to implementation aspects in order to speed up the computation of the sensitivity.

\subsection{Discretization} \label{subsec:discr}

We base our discretization choices on \cite{panzeri2023lab,Fumagalli2023} and here we recall the main aspects. For a more detailed and complete
discussion refer to the aforementioned works.

For the approximation of the domain $\Upsilon$, the liner $\lambda$ and the electrodes $\gamma$, we consider the computational grids made of simplices build with Gmsh \cite{Geuzaine2009}. For the numerical approximation of the problem, for each dimension, we utilize the multi-point
flux approximation (MPFA) scheme \cite{Aavatsmark2007, Aavatsmark2002}, which is a finite volume scheme where the degrees of freedoms are one for each cell of the grid. MPFA solves the problem in primal form, is convergent and guarantees local conservation of the current density for each grid cell, as discussed in the
aforementioned works,
which is an important property for our application, in particular in the computation of the sensitivity. To compute the local gradient of the solution, we first calculate the normal component of the current density at each face and then interpolate over the lowest-order Raviart-Thomas space \cite{Raviart1977}.

\begin{figure}
    \centering
    \includegraphics[width=12cm]{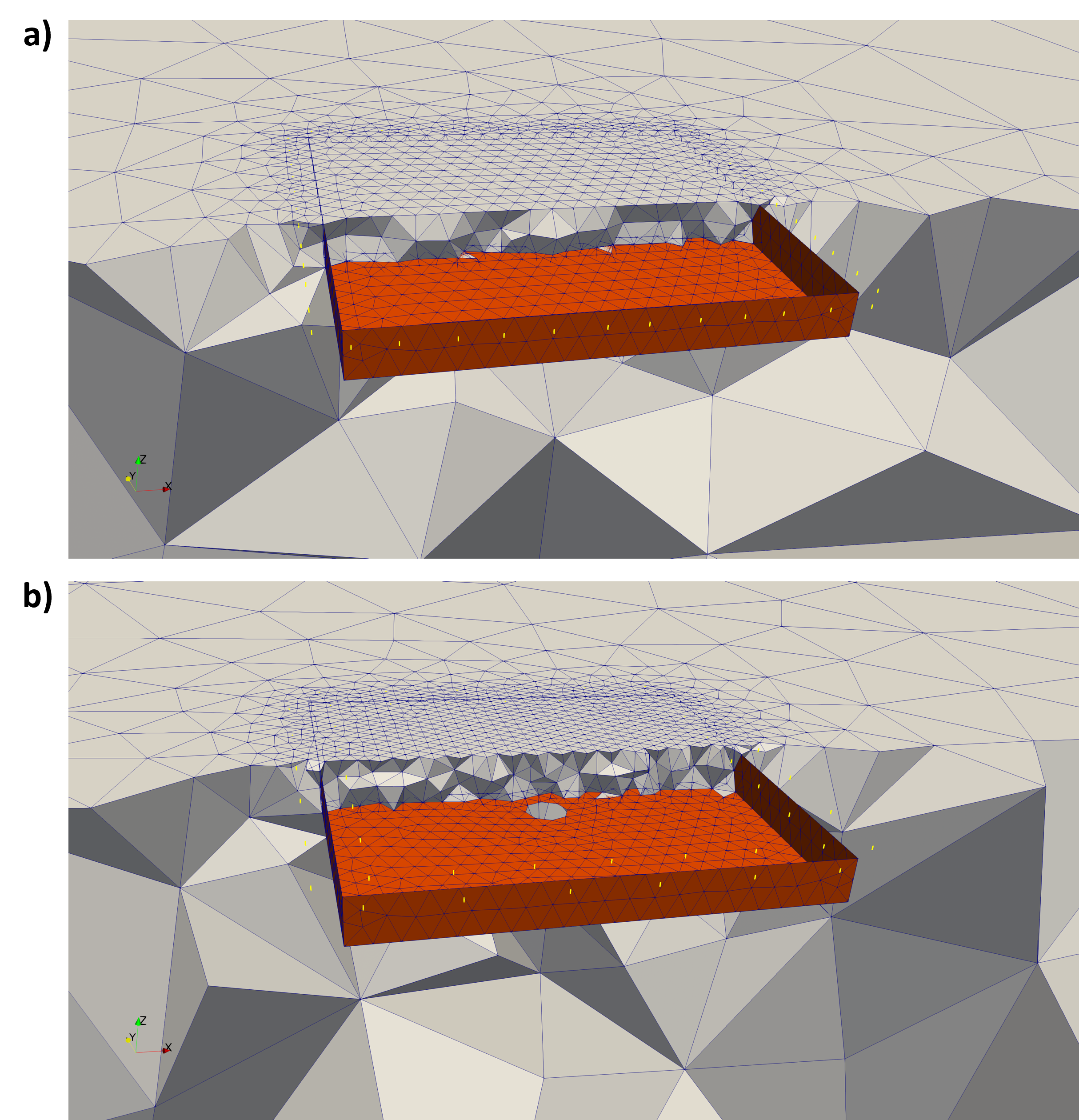}
    \caption{Examples of computational grids for the mixed dimensional model: 3D tetrahedral elements in grey, 2D triangular elements of the liner in red and 1D electrodes in yellow. a) 48 electrodes around the box-shaped liner; b) 24 electrodes outside the liner and 24 electrodes inside, a hole in the bottom surface of the liner with its vertical walls not in contact with the top boundary of the domain. The liner has dimensions $1\sib{\meter}\times1\sib{\meter}\times0.1\sib{\meter}$.}
    \label{fig:grid_zoom2}
\end{figure}

A key aspect in the discretization is associated with the coupling between the dimensions, in particular we consider a conforming mortar-based strategy
for the coupling between the liner and the surrounding material. This means that we consider additional grids at each side of the
liner, conforming to the liner grid, where the interface variables ($j_\lambda$ or $w_\lambda$) are defined. Their degrees of freedom
is one per cell and represent the interface variable. The three-dimensional grid is conforming with the liner grid, see an example reported in
Figure \ref{fig:grid_zoom2}.

To achieve coupling between the electrodes and the domain $\Upsilon$, we embed the computational grids of the electrodes into the grid cells of the domain. This is an effective strategy to reduce the number of three-dimensional cells needed for the coupling but requires a dedicated approach for the discretization of the
interface law. As done before, we consider a mortar-based scheme where a single grid is constructed conforming with the electrode grid but embedded in the three-dimensional grid. Also in this case, we define the interface variables ($j_\gamma$) on the mortar grid, with one degree of freedom per cell.

By employing the discretization schemes, we construct a linear system with a saddle-point structure. With an abuse of notation, we indicate the discrete unknowns with the same symbols of their continuous counterparts. The problem now reads: find $(\varphi, \varphi_\lambda, \varphi_\gamma, j_\lambda, j_\gamma)$ solution of the following linear system
\begin{gather}\label{eq:matrix_form}
    \begin{bmatrix}
        A         &           &          & - B_\lambda^\top & -B_\gamma^\top \\
                  & A_\lambda &          & -C_\lambda^\top  &                \\
                  &           & A_\gamma &                  & -C_\gamma^\top \\
        B_\lambda & C_\lambda &          & M_\lambda                         \\
        B_\gamma  &           & C_\gamma &                  & M_\gamma
    \end{bmatrix}
    \begin{bmatrix}
        \varphi         \\
        \varphi_\lambda \\
        \varphi_\gamma  \\
        j_\lambda       \\
        j_\gamma
    \end{bmatrix}
    =
    \begin{bmatrix}
        0 \\ 0 \\ j_{bc} \\ 0 \\0
    \end{bmatrix}
\end{gather}
where the $A$-matrices represent the stiffness matrices associated to each dimension, the $M$-matrices are mass matrices for the mortar variables, the $B$-matrices are the coupling between the mortars and the three-dimensional variables, and the $C$-matrices are the coupling between the mortars and the lower-dimensional variables. Finally, at the right-hand side, $j_{bc}$ represents the current density injected by the current electrodes, approximation of the boundary condition in \eqref{eq:the_system_electrode}.
The computation of the $\sigma$-weighted Green functions in \eqref{eq:greenmd} requires the solution of a linear system where the matrix is the same as in \eqref{eq:matrix_form} and the right-hand side term is associated to a specific potential electrode.

\subsection{Implementation aspects} \label{subsec:imple}

The computation of the sensitivity, either for single or multiple electrodes and both for the equi-dimensional and the mixed-dimensional cases (see \eqref{eq:sensi}, \eqref{eq:quadrupole}, \eqref{eq:sens_mixed_out} and \eqref{eq:sens_mixed_int}), requires the evaluation of the Green's functions associated to the considered electrodes. Moreover, it is worth pointing out that the classical direct current problem \eqref{eq:the_system} and the sensitivity analysis \eqref{eq:sens} have very similar equations and can be both expressed in terms of
the Green's functions.

By exploiting this similarity, in our code we compute the Green's functions for all the deployed electrodes and then we simply linearly combine them to compute the voltage diﬀerence as well as the sensitivity for each considered electrode conﬁguration. As mentioned above, for several configurations, the global sensitivity is obtained by summing the absolute values of the sensitivity for each configuration.


Another way to speed up the code is to observe that the Green's functions all solve the same problem, similar to the one in \eqref{eq:matrix_form}, but with different source term positions. This means that, thanks to the chosen numerical scheme, we can assemble the discretization matrix only once and select a matrix as the right-hand side, where each column represents a unit source term associated with an electrode. By considering standard software packages, e.g. SuperLU \cite{Li2005}, we can speed up the solution of linear problem.

Finally, at a practical level, imposing the condition \eqref{eq:at_inf} may not be feasible. Therefore, we limit the computational domain and instead impose a homogeneous boundary condition on $\varphi$ to approximate the real condition \eqref{eq:at_inf}. This approximation becomes more accurate as the computational domain increases, but at the cost of greater computational effort. Since the accuracy of the electric potential is less critical far from the electrodes and liner, we coarsen the grid to achieve a more manageable computational cost.

\section{Sensitivity synthetic tests}

The mathematical approach regarding the direct current problem described in the previous sections has been validated and tested by \cite{Fumagalli2023}.
The validation process involved comparing the solver against established analytical solutions (\cite{Sheriff1990}; \cite{aldridge1989direct}), other software (\cite{RUCKER2017106}), as well as laboratory experiments (\cite{Fumagalli2023}; \cite{panzeri2023lab}). Testing, even in very simple setups, confirmed the ability of geoelectrical surveys to detect damaged liners, paving the way for more complex investigations in real-world applications.

Following a similar approach, here we first validate our 3D code for geoelectrical sensitivity analysis by considering a few common quadrupoles. Next, we conduct numerical tests in downscaled and simplified scenarios that are representative of geoelectrical investigations of landfills. Downscaling the domain does not impact computational costs because the grid cells and electrode spacing are proportionally reduced; this step is primarily intended for future comparison with laboratory tests. The domain size is carefully chosen to be large enough to avoid any boundary effects in all simulations. Accordingly, the depth is set to two times the spacing between the two furthest electrodes.

In all tests, the domain is discretized with an unstructured grid with the software Gmsh \cite{Geuzaine2009} and the modelled sensitivity values are normalized by the volume of each cell. The injected current is $1\sib{\ampere}$.

\subsection{4-electrode arrays}

We consider the Wenner-alpha and dipole-dipole (n = 1) arrays because they are the most commonly used quadrupoles due to their versatility and effectiveness. Additionally, and most interestingly, the Wenner array has a strong ability to detect horizontal structures, as its sensitivity exhibits a marked gradient with depth. In contrast, the dipole-dipole array is most sensitive to changes in resistivity along the horizontal direction.

For both arrays, electrodes are placed on the surface of a 3D homogeneous half-space with dimensions of $6\sib{\meter}\times4\sib{\meter}\times2\sib{\meter}$. Assuming the origin of the reference system is at the lower bottom corner of the domain, the coordinates of the furthest electrodes are (2, 2, 2) $\sib{\meter}$ and (4, 2, 2) $\sib{\meter}$, respectively, with an inter-electrode distance of $0.66\sib{\meter}$. Electrode configurations are C1-P1-P2-C2 and C2-C1-P1-P2 for Wenner-alpha and dipole-dipole arrays, respectively, being C1 and P1 positive current and potential electrodes. To obtain accurate results without excessively increasing computational costs, we set the characteristic length of the tetrahedral mesh elements to $0.05\sib{\meter}$ near the electrodes and $0.5\sib{\meter}$ at the boundaries of the domain. For validation purposes, we consider the analytical expression of the pole-pole sensitivity \eqref{eq:sensi} described in \cite{loke1995least}.

\begin{figure}[H]
    \centering
    \includegraphics[width=15cm]{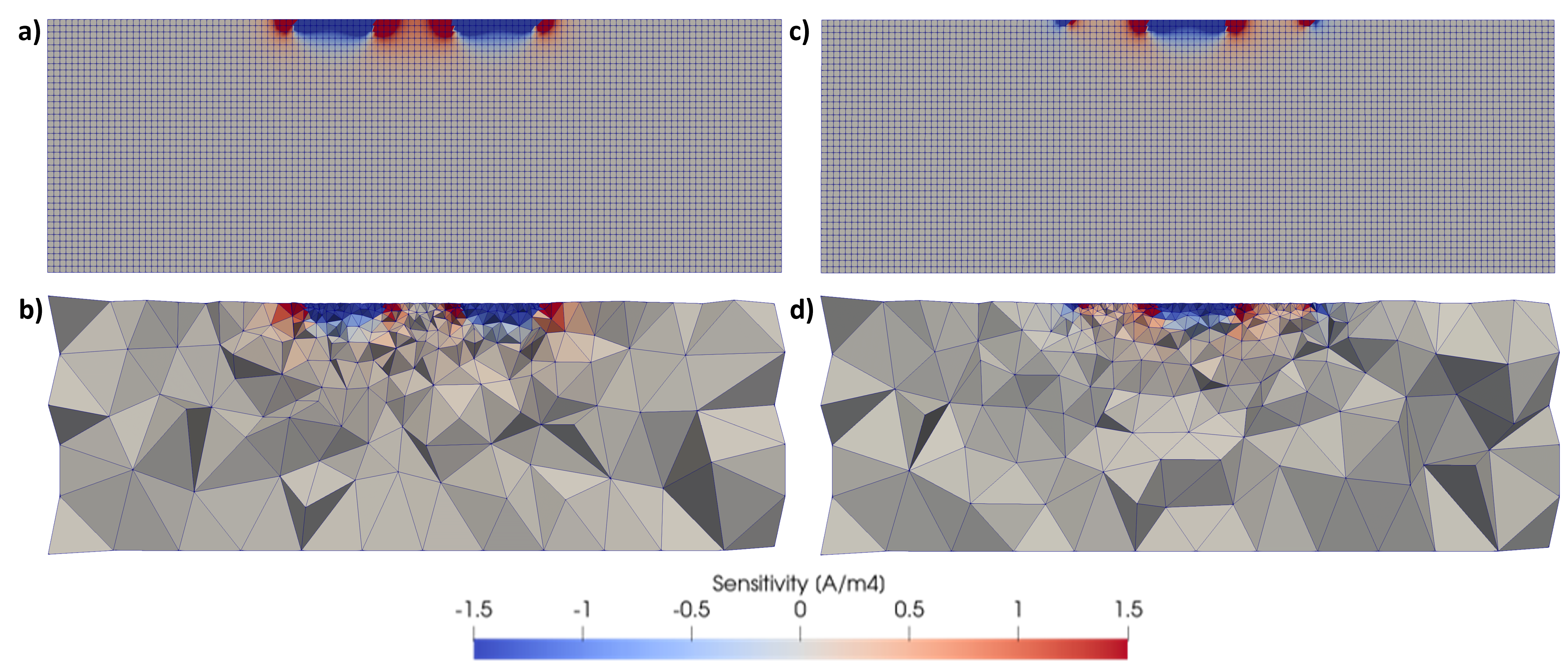}
    \caption{Vertical slices of sensitivity values for Wenner-alpha a)-b) and dipole-dipole c)-d) arrays and a homogeneous half-space. a)-c) show plots of the analytical expression \eqref{eq:sensi}, while b)-d) show the modelled results. Both arrays have total length of $2\sib{\meter}$.}
    \label{fig:analytical_modelled_WA_dip_dip}
\end{figure}

Figure \ref{fig:analytical_modelled_WA_dip_dip} shows the sensitivity values on vertical slices extracted along the tested arrays from the 3D domain. Figures \ref{fig:analytical_modelled_WA_dip_dip}a-c are obtained using the analytical expression evaluated on a structured grid with $0.05\sib{\meter}$ spatial sampling, while Figures \ref{fig:analytical_modelled_WA_dip_dip}b-d show the modeled results. It is worth noting that the structured grids contain about 340000 elements, compared to only 9500 elements in the unstructured grids.

\subsection{48-electrode arrays and a simplified landfill}

For our downscaled tests, we set a 3D domain with dimensions $10\sib{\meter}\times10\sib{\meter}\times3\sib{\meter}$, with a $1\sib{\meter}\times1\sib{\meter}\times0.1\sib{\meter}$ box-shaped liner placed at the top center (i.e., the vertical walls of the box are in contact with the top surface of the domain; see Figure \ref{fig:grid_zoom2}a). Obviously, the liner is modelled as a 2D object by our mixed-dimensional code. The resistivity values inside and outside the liner (i.e., the landfill) are set to $20\sib{\ohm\meter}$ and $100 \sib{\ohm\meter}$, respectively \cite{bernstone2000dc}. The resistivity of the liner is set to $10^{15} \sib{\ohm\meter}$ in agreement with the electrical properties of HDPE.
The characteristic length of the computational grid is set to $0.02\sib{\meter}$ near the electrodes, to $0.05\sib{\meter}$ for the liner and to $2\sib{\meter}$ at the boundaries of the domain (Figure \ref{fig:grid_zoom2}).
We assume to use 48 channels, as this is standard for most multi-electrode resistivity meters on the market, deployed at the surface. In general, electrodes placed near the perimeter edges of the liner are preferable to those positioned across the landfill body, as the latter, because of landfill elevation, are farther from the bottom of the landfill, resulting in reduced penetration and lower resolution at depth.

Three different settings are defined: \textit{case 1)} 48 electrodes deployed along the outer perimeter of the liner with a spacing of $0.08 \sib{\meter}$ (Figure \ref{fig:grid_zoom2}a); \textit{case 2)} 24 electrodes with $0.16 \sib{\meter}$ spacing along the outer perimeter of the liner and 24 electrodes with $0.12 \sib{\meter}$ spacing along the inner perimeter of the liner (Figure \ref{fig:grid_zoom2}b); \textit{case 3)} 24 electrodes with $0.16 \sib{\meter}$ spacing along the outer perimeter of the liner and 24 electrodes along a grid inside the liner ($0.2 \sib{\meter}$ spacing along x and $0.14 \sib{\meter}$ spacing along y).
All cases involve about 15000 cells and more than 635000 configurations for the evaluation of the sensitivitiy, considering all possible electrode combinations and permutations
(excluding reciprocal configurations) with geometric factors smaller than $10^4\sib{\meter}$, including pole-dipole and pole-pole arrays. For all the settings, we also considered a test with a $0.1\sib{\meter}$ diameter hole in the bottom surface of the liner (Figure \ref{fig:grid_zoom2}b).

\begin{figure}[H]
    \centering
    \includegraphics[width=16cm]{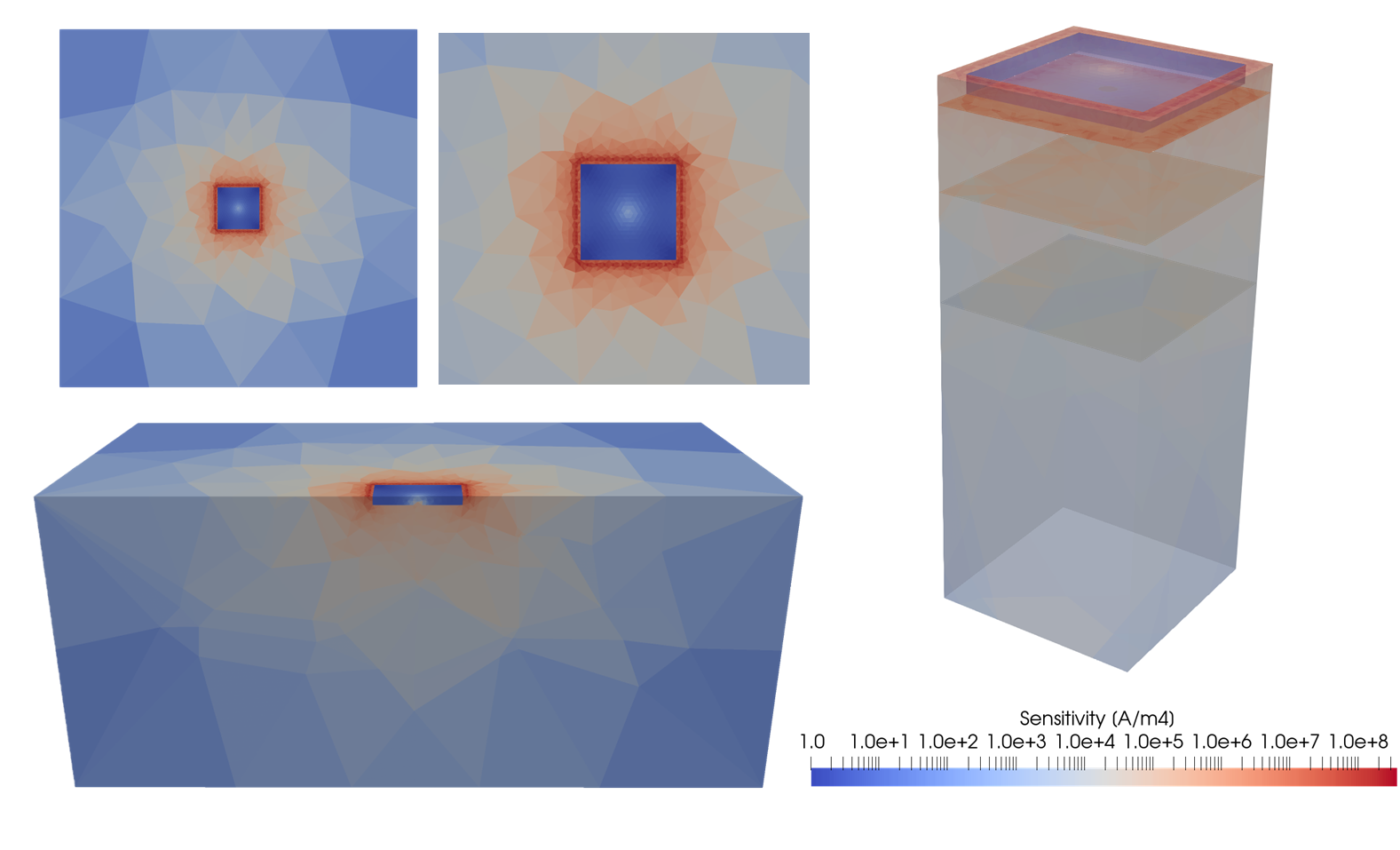}
    \caption{Images of computed sensitivity for \textit{case 1)} with a $0.1\sib{\meter}$ diameter hole in the bottom surface of the liner. From top-left, clockwise: top view, zoomed top view, volume below the box-shaped liner with 3 depth slices, vertical cut. See text for details.}
    \label{fig:48_ouside_hole}
\end{figure}

\begin{figure}[H]
    \centering
    \includegraphics[width=15cm]{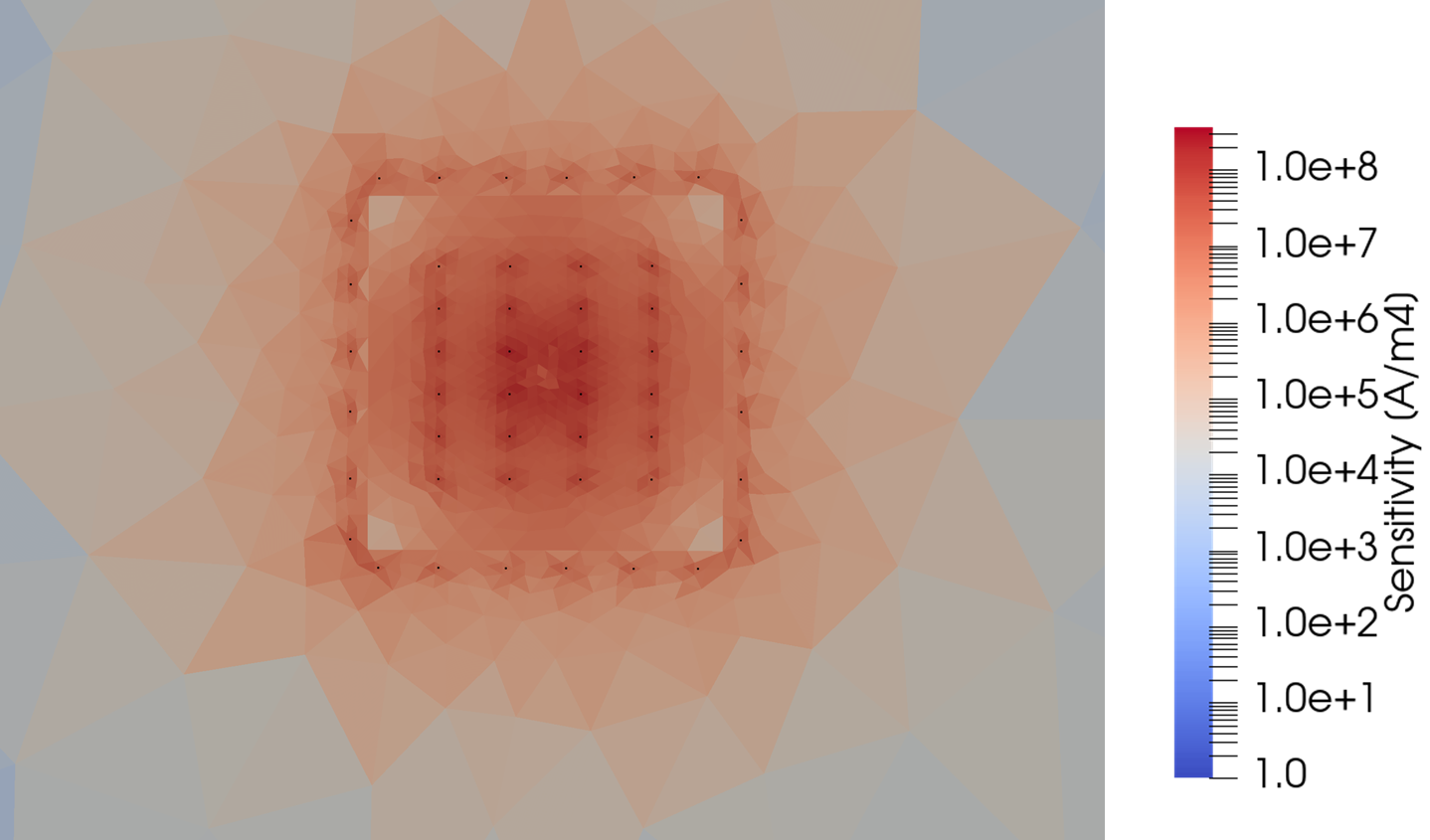}
    \caption{Zoomed top view of the computed sensitivity for \textit{case 3)} with a hole in the middle of the liner.}
    \label{fig:grid_zoom}
\end{figure}

Figures \ref{fig:48_ouside_hole} and \ref{fig:grid_zoom} show some images of the sensitivity computed for \textit{cases 1)} and \textit{3)}.

In addition, \textit{cases 2)} and \textit{3)} are tested with the box-shaped liner shifted downward by $0.05 \sib{\meter}$, so its vertical walls are not in contact with the top boundary of the 3D domain (Figure \ref{fig:grid_zoom2}b). This setup simulates the condition of an electrical connection between the landfill and the surrounding media that is not due to damages to the liner. The electrical connection could be caused by movements of the liner due to the weight of the waste, or by materials that can accidentally cover the perimeter edges of the liner.


\section{Discussion}

\begin{figure}[H]
    \centering
    \includegraphics[width=13cm]{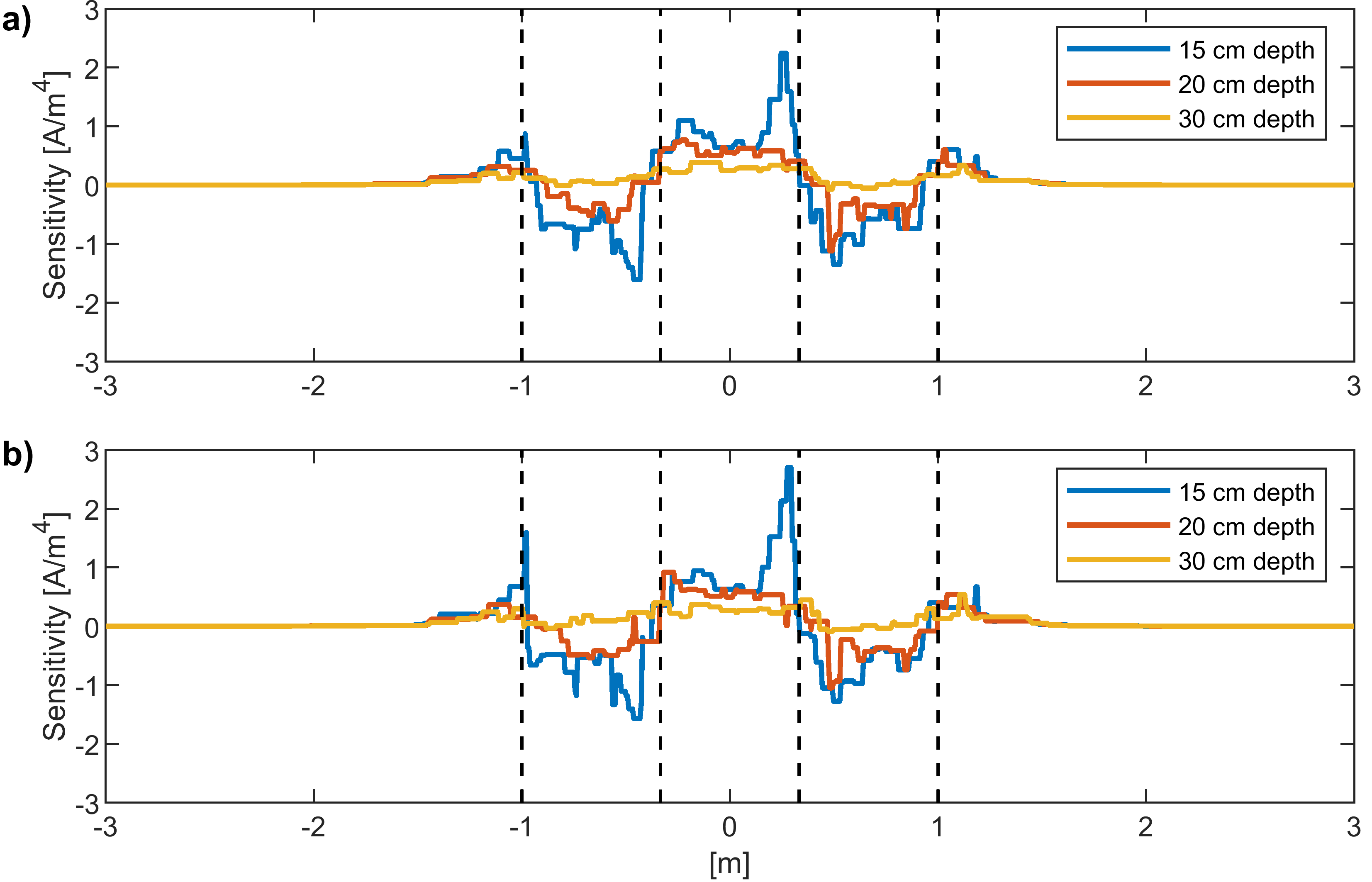}
    \caption{Sensitivity values along horizontal lines at different depths for the Wenner-alpha array. Vertical dashed lines indicate the positions of the electrodes (from right: C1-P1-P2-C2). a) Analytical values; b) modelled values.}
    \label{fig:horizontal plot over line_WA}
\end{figure}

\begin{figure}[H]
    \centering
    \includegraphics[width=13cm]{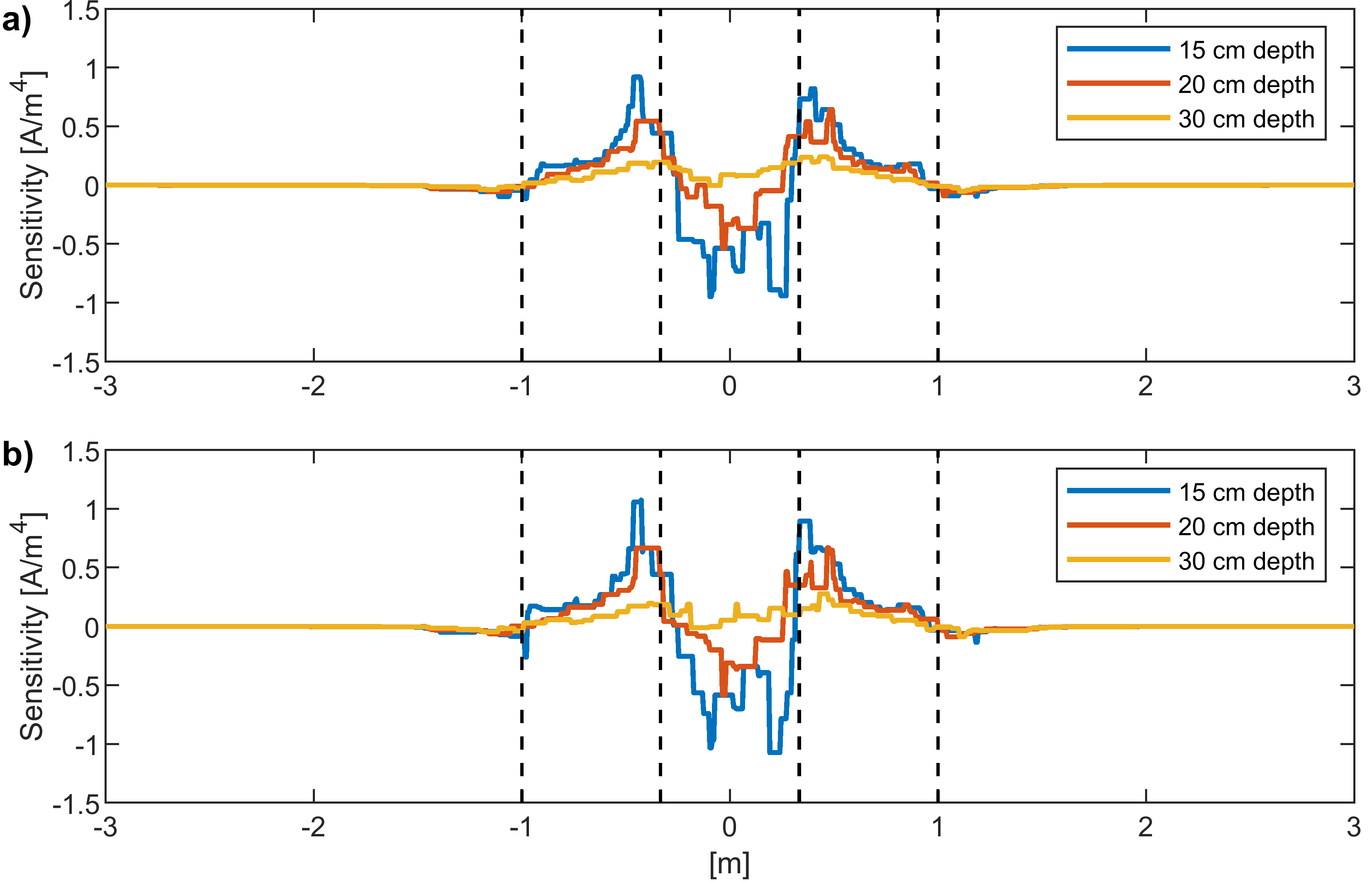}
    \caption{Sensitivity values along horizontal lines at different depths for the dipole-dipole array (n = 1). Vertical dashed lines indicate the positions of the electrodes (from right: C2-C1-P1-P2). a) Analytical values; b) modelled values.}
    \label{fig:horizontal plot over line_DD}
\end{figure}

As far as validation is concerned, Figures \ref{fig:horizontal plot over line_WA} and \ref{fig:horizontal plot over line_DD} report the results obtained for the Wenner-alpha and dipole-dipole arrays, respectively. More in detail, the sensitivity values computed with our code are evaluated along horizontal lines at $0.15 \sib{\meter}$, $0.2 \sib{\meter}$ and $0.3 \sib{\meter}$ depth below the arrays and analytical values are computed at the same positions with \eqref{eq:sensi}, as mentioned in the previous paragraph. The irregular stairstep plots are due to the different sizes of the mesh elements intercepted by the horizontal lines (see Figures \ref{fig:analytical_modelled_WA_dip_dip}b-d).

We can observe a very good match between simulated and analytical values, with percent root-mean-square error (RMSE) around 8 \% and 4 \% at $0.15 \sib{\meter}$ depth for Wenner-alpha and dipole-dipole arrays, respectively. As expected, peaks are recorded near the electrodes at shallow depths. By considering even shallower depths, the amplitude of the peaks increases, and the mismatch between modeled and analytical values tends to grow due to the difficulty of the numerical code in handling the singularity of the solution at a reasonable computational cost. On the contrary, percent RMSE is below 1 \% at a depth of $0.3 \sib{\meter}$. At this depth (approximately half the electrode spacing), a significant decrease in sensitivity is observed across the entire domain, consistent with the median depth of investigation (\cite{Edwards1977}), which is about 15\% of the total length for both arrays (\cite{Loke2022}). The plots also indicate a higher vertical sensitivity gradient for the Wenner-alpha array, related to its high vertical resolution. In contrast, the pronounced lateral sensitivity gradients of the dipole-dipole array confirm its good horizontal resolution.

Regarding the results obtained from the simplified landfill setups, we first note that the outcomes of \textit{case 3)} do not differ significantly from those of \textit{case 2)}; therefore, they will not be discussed further.
Second, we note that computing the global sensitivity as the sum of the absolute values of the sensitivities associated with each configuration is a simple but rough method, as it does not account for the fact that information from two measurements of the same region might not be independent. A more mathematically rigorous (and computationally intensive) approach would involve, for instance, considering the model resolution matrix, which relates the estimated resistivity values to the actual values in the geoelectrical inverse problem (\cite{Day-Lewis2005}).

\begin{figure}[H]
    \centering
    \includegraphics[width=15cm]{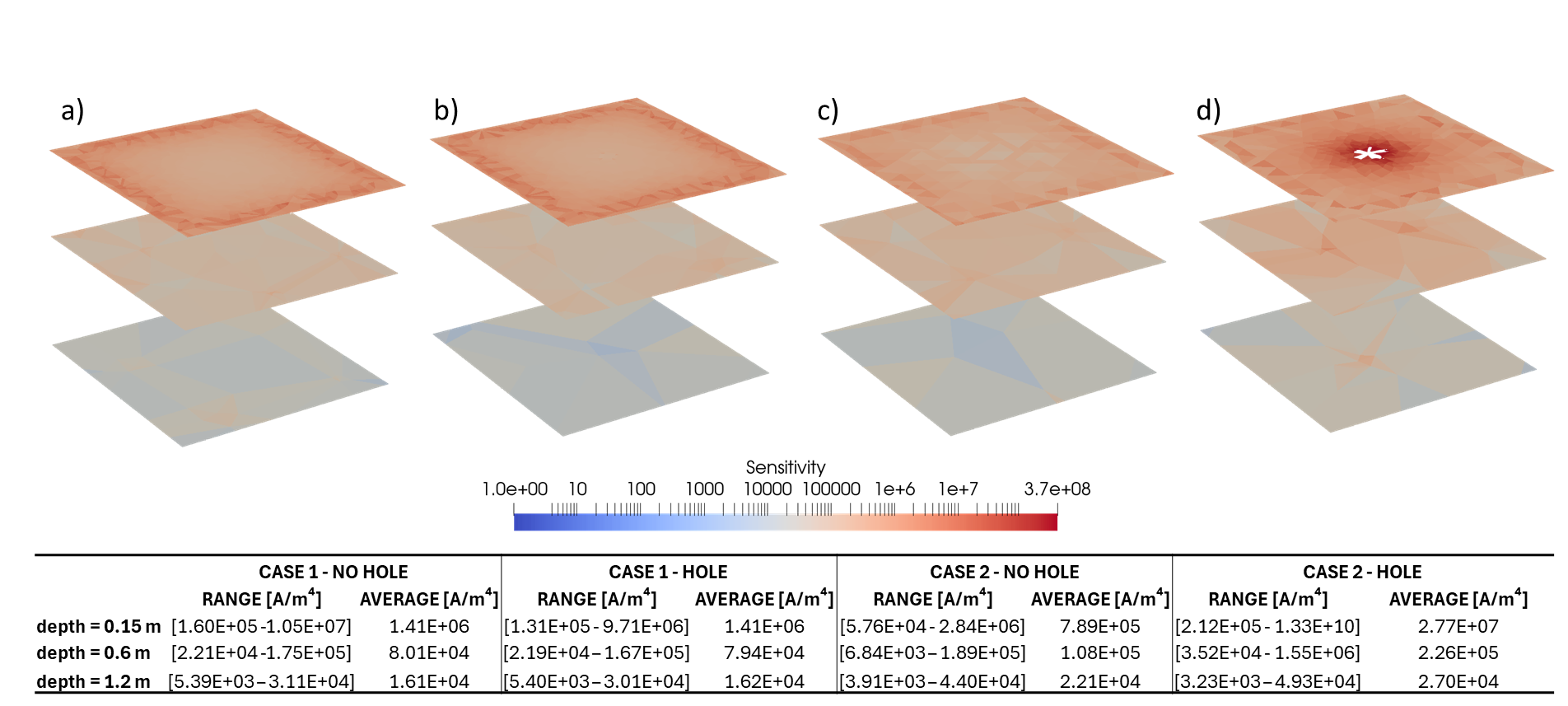}
    \caption{Sensitivity slices at $0.15 \sib{\meter}$, $0.6 \sib{\meter}$ and $1.2 \sib{\meter}$ depth for different setups. a) 48 electrodes outside the box-shaped liner (\textit{case 1)}); b) same as a) but with a hole in the liner; c) 24 electrodes outside the box-shaped liner and 24 electrodes inside (\textit{case 2)}); d) same as c) but with a $0.1\sib{\meter}$ diameter hole in the bottom surface of the liner.}
    \label{fig:integratevariable}
\end{figure}

As expected, the images show that sensitivity is highest near the electrodes and is very low inside the box-shaped liner when no electrodes are present due to the high resistivity of HDPE (Figures \ref{fig:48_ouside_hole} and \ref{fig:grid_zoom}).
To evaluate the sensitivity values in the most relevant area, we extracted three $1\sib{\meter}\times1\sib{\meter}$ horizontal slices at depths of $0.15 \sib{\meter}$, $0.6 \sib{\meter}$ and $1.2 \sib{\meter}$ below the liner for all setups. These slices correspond to the area immediately below the liner and to depths roughly equal to half and one distance between the furthest electrodes. Figure \ref{fig:integratevariable} presents the slices along with a table that lists, for each slice, the range of volume-normalized sensitivity values and the area-weighted average sensitivity.
In the case of 48 electrodes placed outside the liner (\textit{case 1)}), the results do not significantly change if the liner is damaged, as the current is forced to flow between the electrodes without passing through the high-resistivity, box-shaped structure. On the contrary, for the case with 24 electrodes inside and 24 electrodes outside (\textit{case 2)}), higher values are recorded at all depths when there is a hole, especially immediately below the liner, where the average and maximum sensitivity values increase of 2 and 4 orders of magnitude, respectively.
It is interesting to note that, if there is no hole, sensitivity values of the shallowest slice are better in \textit{case 1)}. This is reasonably due to the higher number of electrodes outside the box-shaped liner that forces current flow in that region.

\begin{figure}[H]
    \centering
    \includegraphics[width=12cm]{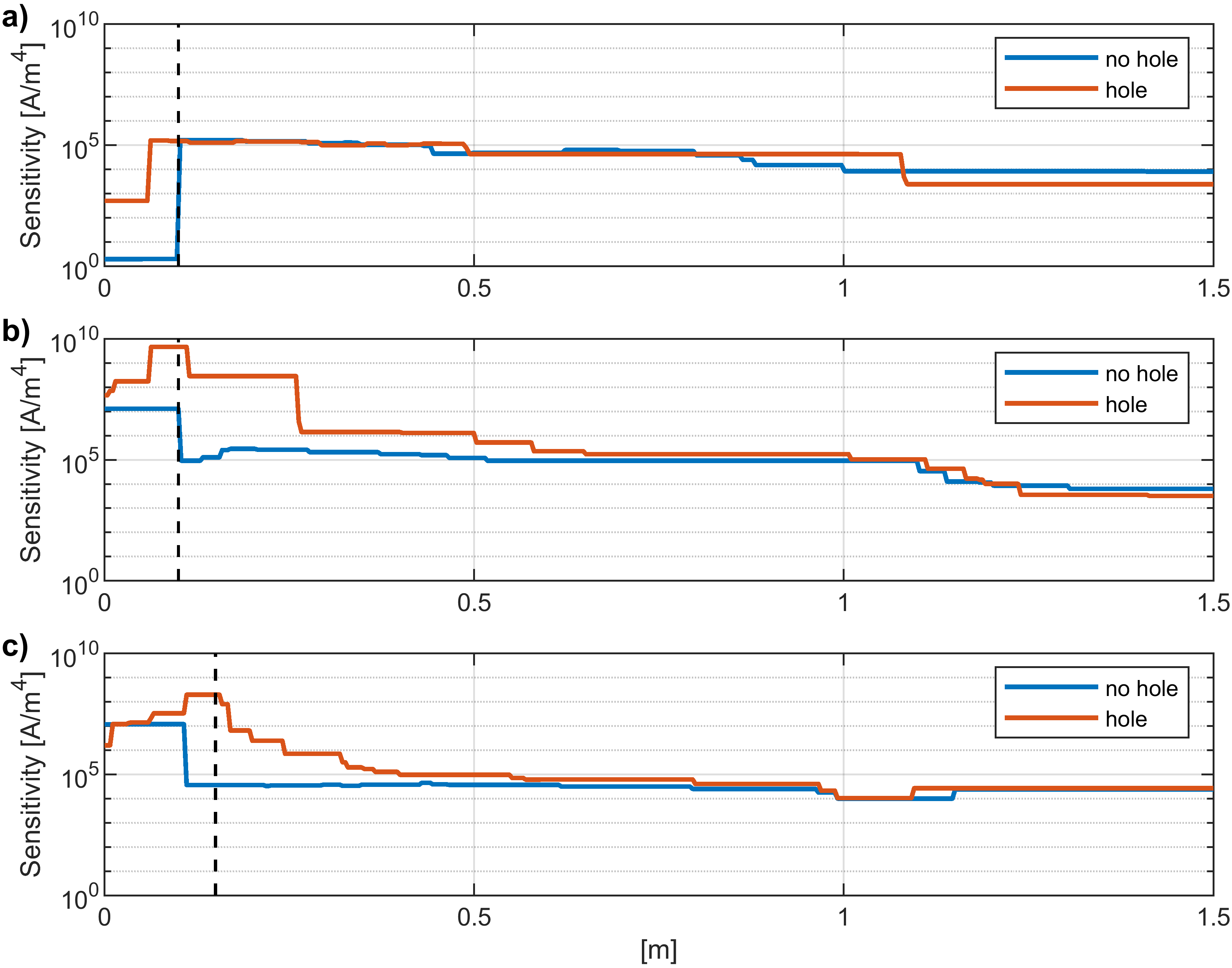}
    \caption{Sensitivity values along a vertical
line at the center of the domain. a) 48 electrodes outside the box-shaped liner (\textit{case 1)}); b) 24 electrodes outside the box-shaped liner and 24 electrodes inside (\textit{case 2)}); c) same as b) but with liner shifted down by $0.05\sib{\meter}$. Vertical dashed line indicates the depth of the bottom surface of the liner.}
    \label{fig:plotoverlinevertical}
\end{figure}

Similar considerations can be made by examining the trend of sensitivity values along a vertical line at the center of the domain (Figure \ref{fig:plotoverlinevertical}). The line goes from the top to the bottom of the domain and  passes through the hole when the hole is present.
When there is no hole, sensitivity values are close to zero inside the box-shaped liner only for \textit{case 1)} (Figure \ref{fig:plotoverlinevertical}a), while are similar for all tested cases below the liner.
Conversely, the tested cases show significant differences when the hole is present. In \textit{case 1)}, sensitivity values increase just inside the liner, whereas in \textit{case 2)}, there is a marked increase in sensitivity across and immediately below the liner (Figure \ref{fig:plotoverlinevertical}b). This trend is also observed when considering the setting with the liner shifted down by $0.05\sib{\meter}$ (Figure \ref{fig:plotoverlinevertical}c).
From these results, the arrangement of electrodes as in \textit{case 2)} may be promising for detecting possible liner damage, even if there is an electrical connection between the landfill and the surrounding media due to improper membrane deployment.

Additional tests performed with smaller holes, with diameters down to 0.02 m, do not provide any new insights. Any hole in the liner can cause a significant increase in sensitivity values compared to the case without a hole, particularly near the liner. This occurs because significant current densities are concentrated in the area of the hole, as shown by \cite{Fumagalli2023}, and sensitivity increases near any current source (see, for example, Figure \ref{fig:analytical_modelled_WA_dip_dip}).

\section{Conclusions}

We have presented and validated a new mixed-dimensional code to accurately compute sensitivity of direct current investigations of MSWLFs, where a highly resistive membrane is typically placed underneath the waste mass to avoid subsurface contamination concerns. The mixed-dimensional framework is useful to tackle the computational costs and ill-posedness of the geoelectrical inverse problem by modelling the electrodes and, most importantly, the resistive liner, as 1D and 2D elements, respectively.
Besides some technical implementation aspects described before, the code uses Green's functions to speed up computation of both voltage difference and the sensitivity for each considered electrode configuration.

According to preliminary modeling of downscaled and simplified landfill settings, acquisitions with all electrodes placed outside the perimeter edges of the box-shaped liner may not be able to detect a hole in the most challenging location, i.e., at the center of the liner’s bottom surface. Conversely, acquisitions with the same number of electrodes deployed along the perimeter of the liner, but placed both inside and outside the landfill, may detect potential damage to the liner as sensitivity increases near the damaged area. No significant improvement has been observed when the electrodes inside the landfill are arranged in a grid pattern. Therefore, electrodes arranged linearly on both sides of the liner should be preferred because their deployment is logistically easier, particularly in the presence of steep and uneven landfill topography. In addition, electrodes running across the landfill topography may be too far from the underlying liner, leading to issues with resolution and penetration depth. Furthermore, the topography may affect the final results if not properly accounted for. Electrodes arranged linearly along both sides of the perimeter edges of the liner can also be used to check the electrical insulation of the landfill relative to the surrounding media by measuring the electrical resistance between adjacent electrode pairs.

Our simulations show that the results do not differ significantly for holes with diameters ranging from approximately one electrode spacing down to one-sixth of the electrode spacing. Even with relatively small damage, the sensitivity increases by 2-3 orders of magnitude, which is promising for landfill monitoring purposes.

The above considerations still apply when the liner is not in contact with the top of the domain. This setting is common in the field and can simulate electrical connections between the landfill and the surrounding media due to downshifts of the liner or conductive materials covering its perimeter edges.

Although the preliminary results obtained here are promising for the geoelectrical monitoring of landfills, it is important to note that we have considered an extremely high number of electrode configurations for the calculation of sensitivity. Moving forward, an appropriate number of measurements will need to be evaluated to ensure sensitivity levels that meet the survey objectives with a reasonable field effort. Our code to compute the sensitivity can be obviously used to select the best configuration and to design electrical acquisitions with the aim of both achieving the appropriate investigation depth (e.g., where the liner is located) and spatial resolution.

Downscaled laboratory tests will be necessary to validate the modeling results and determine whether the computed sensitivity values are sufficiently high to detect liner damage.

To improve sensitivity at depth, at the expense of higher costs and more complex data processing, cross-borehole ERT imaging can also be considered.

Finally, although we have noted that sensitivity affects the resolution and penetration depth of geoelectrical surveys and depends on factors such as physical approximations in the forward model, survey geometry, and resistivity distribution, real scenarios also require consideration of measurement errors, data signal-to-noise ratio, as well as parameterization and regularization used in the inversion (\cite{Binley2005}).

\section*{Acknowledgement}

We wish to thank Prof. Luigi Zanzi for fruitful discussions. The present research is part of the activities of ''Dipartimento di Eccellenza 2023-2027'', Italian Minister of University and Research (MUR), grant Dipartimento di Eccellenza 2023-2027.


\begin{thebibliography}{BDOH00}

\bibitem[Aav02]{Aavatsmark2002}
Ivar Aavatsmark.
\newblock An introduction to multipoint flux approximations for quadrilateral grids.
\newblock {\em Computational Geosciences}, 6:405--432, 2002.

\bibitem[Aav07]{Aavatsmark2007}
Ivar Aavatsmark.
\newblock Multipoint flux approximation methods for quadrilateral grids.
\newblock In {\em The 9\textsuperscript{th} International Forum on Reservoir Simulation, Abu Dhabi}, pages 9--13, 2007.

\bibitem[AHZA20]{Aguzzoli2020}
A.~Aguzzoli, A.~Hojat, L.~Zanzi, and D.~Arosio.
\newblock Two dimensional ert simulations to check the integrity of geomembranes at the base of landfill bodies.
\newblock 2020(1):1--5, 2020.

\bibitem[AO89]{aldridge1989direct}
DF~Aldridge and DW~Oldenburg.
\newblock Direct current electric potential field associated with two spherical conductors in a whole-space 1.
\newblock {\em Geophysical prospecting}, 37(3):311--330, 1989.

\bibitem[BD03]{Binley2003}
Andrew Binley and William Daily.
\newblock The performance of electrical methods for assessing the integrity of geomembrane liners in landfill caps and waste storage ponds.
\newblock {\em Journal of Environmental and Engineering Geophysics}, 8(4):227--237, 2003.

\bibitem[BDOH00]{bernstone2000dc}
Christian Bernstone, Torleif Dahlin, Tomas Ohlsson, and H~Hogland.
\newblock Dc-resistivity mapping of internal landfill structures: two pre-excavation surveys.
\newblock {\em Environmental Geology}, 39:360--371, 2000.

\bibitem[BK05]{Binley2005}
Andrew Binley and Andreas Kemna.
\newblock {\em DC Resistivity and Induced Polarization Methods}, pages 129--156.
\newblock Springer Netherlands, Dordrecht, 2005.

\bibitem[Com15]{EC2015}
European Commission.
\newblock Closing the loop—an eu action plan for the circular economy. communication from the commission to the european parliament, the council, the european economic and social committee and the committee of the regions.
\newblock Technical report, Brussels: European Commission., 2015.

\bibitem[DC17]{DEDONNO2017302}
Giorgio {De Donno} and Ettore Cardarelli.
\newblock Tomographic inversion of time-domain resistivity and chargeability data for the investigation of landfills using a priori information.
\newblock {\em Waste Management}, 59:302--315, 2017.

\bibitem[DLB05]{Day-Lewis2005}
K.~Singha Day-Lewis, F.~D. and A.~M. Binley.
\newblock Applying petrophysical models to radar travel time and electricalresistivity tomograms: Resolution-dependent limitations.
\newblock {\em J. Geophys. Res.}, 110:B08206, 2005.

\bibitem[DPC{\etalchar{+}}13]{DeCarlo2013}
Lorenzo {De Carlo}, Maria~Teresa Perri, Maria~Clementina Caputo, Rita Deiana, Michele Vurro, and Giorgio Cassiani.
\newblock Characterization of a dismissed landfill via electrical resistivity tomography and mise-{\'a}-la-masse method.
\newblock {\em Journal of Applied Geophysics}, 98:1--10, 2013.

\bibitem[DRL10]{Dahlin2010}
T.~Dahlin, H.~Rosqvist, and V.~Leroux.
\newblock Resistivity-ip mapping for landfill applications.
\newblock {\em First Break}, 28(8), 2010.

\bibitem[ea13]{Jones2013}
Peter Tom~Jones et~al.
\newblock Enhanced landfill mining in view of multiple resource recovery: a critical review.
\newblock {\em Journal of Cleaner Production}, 55:45--55, 2013.

\bibitem[FPF{\etalchar{+}}23]{Fumagalli2023}
Alessio Fumagalli, Lorenzo Panzeri, Luca Formaggia, Anna Scotti, and Diego Arosio.
\newblock A mixed-dimensional model for direct current simulations in the presence of a thin high-resistivity liner.
\newblock {\em International Journal for Numerical Methods in Engineering}, 2023.
\newblock Cited by: 0; All Open Access, Green Open Access.

\bibitem[Fra97]{Frangos1997}
William Frangos.
\newblock Electrical detection of leaks in lined waste disposal ponds.
\newblock {\em Geophysics}, 62(6):1737--1744, 1997.

\bibitem[Gir14]{Giroux2014}
Laurie Giroux.
\newblock State of waste management in canada. giroux environmental consulting.
\newblock Technical report, Canadian Council of Ministers of the Environment. Kanata, Ontario, Canada, 2014.

\bibitem[GR09]{Geuzaine2009}
Christophe Geuzaine and Jean-Fran{\c{c}}ois Remacle.
\newblock Gmsh: A 3-d finite element mesh generator with built-in pre- and post-processing facilities.
\newblock {\em International Journal for Numerical Methods in Engineering}, 79(11):1309--1331, 2009.

\bibitem[LB95]{loke1995least}
MH~Loke and RD~Barker.
\newblock Least-squares deconvolution of apparent resistivity pseudosections.
\newblock {\em Geophysics}, 60(6):1682--1690, 1995.

\bibitem[Li05]{Li2005}
Xiaoye~S. Li.
\newblock An overview of {SuperLU}: Algorithms, implementation, and user interface.
\newblock {\em Transactions on Mathematical Software}, 31(3):302--325, September 2005.

\bibitem[Lok22]{Loke2022}
Meng~Heng Loke.
\newblock Tutorial : 2-d and 3-d electrical imaging surveys.
\newblock 2022.

\bibitem[LRQ{\etalchar{+}}19]{Ling2019}
C.~Ling, A.~Revil, Y.~Qi, F.~Abdulsamad, P.~Shi, S.~Nicaise, and L.~Peyras.
\newblock Application of the mise-{\`a}-la-masse method to detect the bottom leakage of water reservoirs.
\newblock {\em Engineering Geology}, 261:105272, Nov 2019.

\bibitem[L.S77]{Edwards1977}
Edwards L.S.
\newblock A modified pseudosection for resistivity and induced-polarization.
\newblock {\em Geophysics}, 42(5):939--1087, 1977.

\bibitem[MO90]{mcgillivray1990methods}
Peter~R Mcgillivray and DW~Oldenburg.
\newblock Methods for calculating fr{\'e}chet derivatives and sensitivities for the non-linear inverse problem: A comparative study 1.
\newblock {\em Geophysical prospecting}, 38(5):499--524, 1990.

\bibitem[NB21]{Nanda2021}
Sonil Nanda and Franco Berruti.
\newblock Municipal solid waste management and landfilling technologies: a review.
\newblock {\em Environmental Chemistry Letters}, 19:1433--1456, 2021.

\bibitem[Pet16]{Peters2016}
Adele Peters.
\newblock These maps show how many landfills there are in the u.s., 2016.
\newblock (accessed 20 July 2024).

\bibitem[PFA{\etalchar{+}}23]{panzeri2023lab}
L~Panzeri, A~Fumagalli, A~Aguzzoli, L~Zanzi, L~Longoni, M~Papini, and D~Arosio.
\newblock Lab and modelling dc resistivity tests to analyse the response of a high resistivity liner.
\newblock In {\em 5th Asia Pacific Meeting on Near Surface Geoscience \& Engineering}, volume 2023, pages 1--5. European Association of Geoscientists \& Engineers, 2023.

\bibitem[PFZ{\etalchar{+}}23]{panzeri2023validation}
L~Panzeri, A~Fumagalli, L~Zanzi, L~Longoni, M~Papini, and D~Arosio.
\newblock Validation of a mixed-dimensional code for the analysis of highly resistive liners in landfills.
\newblock In {\em NSG2023 29th European Meeting of Environmental and Engineering Geophysics}, volume 2023, pages 1--5. European Association of Geoscientists \& Engineers, 2023.

\bibitem[PV91]{park1991inversion}
Stephen~K Park and Gregory~P Van.
\newblock Inversion of pole-pole data for 3-d resistivity structure beneath arrays of electrodes.
\newblock {\em Geophysics}, 56(7):951--960, 1991.

\bibitem[RGW17]{RUCKER2017106}
Carsten R{\"u}cker, Thomas G{\"u}nther, and Florian~M. Wagner.
\newblock pygimli: An open-source library for modelling and inversion in geophysics.
\newblock {\em Computers \& Geosciences}, 109:106--123, 2017.

\bibitem[RT77]{Raviart1977}
Pierre-Arnaud Raviart and Jean-Marie Thomas.
\newblock A mixed finite element method for second order elliptic problems.
\newblock {\em Lecture Notes in Mathematics}, 606:292--315, 1977.

\bibitem[Sal16]{Salsa2016}
Sandro Salsa.
\newblock {\em Partial Differential Equations in Action. From Modelling to Theory}, volume~99 of {\em La Matematica per il 3+2}.
\newblock Springer International Publishing, 2016.

\bibitem[TVFT14]{Tsourlos2014}
P.~Tsourlos, G.N. Vargemezis, I.~Fikos, and G.N. Tsokas.
\newblock Dc geoelectrical methods applied to landfill investigation: case studies from greece.
\newblock {\em First Break}, 32(8), 2014.

\bibitem[WR15]{Wagner2015}
Travis~P. Wagner and Tom Raymond.
\newblock Landfill mining: Case study of a successful metals recovery project.
\newblock {\em Waste Management}, 45:448--457, 2015.

\bibitem[WT90]{Sheriff1990}
R.E.~Sheriff W.M.~Telford, L.P.~Geldart.
\newblock {\em Applied Geophysics}.
\newblock Cambridge University Press, 2nd edition, 1990.

\end{thebibliography}

\newcommand{\etalchar}[1]{$^{#1}$}

\printindex

\end{document}